\documentclass[10pt,a4paper,twosides]{article}

\usepackage{amsfonts}
\usepackage{amssymb}
\usepackage{latexsym}
\usepackage{textcomp}
\usepackage[all]{xy}
\usepackage{graphics}
\usepackage{graphicx}

\begin{document}

\newtheorem{theorem}{Theorem}
\newtheorem{observation}{Observation}
\newtheorem{lemma}{Lemma}
\newtheorem{proposition}{Proposition}
\newtheorem{example}{Example}
\newtheorem{corollary}{Corollary}
\newtheorem{definition}{Definition}


%
\title{Topological Classification of Holomorphic, Semi-Hyperbolic Germs, in ``Quasi-Absence'' of Resonances}
\author{Pietro Di Giuseppe \\  \small\textsf{p.digiuseppe@sns.it}}
\normalsize
\date{}
%
%
%
%
\maketitle
\begin{abstract}
The classification, by topological conjugacy, of invertible holomorphic germs $f:(\mathbb{C}^n,0)\rightarrow (\mathbb{C}^n,0)$, with $\lambda_1,\ldots,\lambda_n$ eigenvalues of $df_0$, and $|\lambda_i|\neq 1$ for $i=2,\ldots,n$ while $\lambda_1$ is a root of the unity, is given, in the suitable hypothesis of ``quasi-absence'' of resonances.
\end{abstract}


\section{Introduction}
\label{intro}

In this note we shall deal with germs of holomorphic diffeomorphisms at $0\in\mathbb{C}^n$, understood as discrete holomorphic local dynamical systems. Two such germs $f,g:(\mathbb{C}^n,0)\rightarrow(\mathbb{C}^n,0)$ are locally topologically ($C^r$-diffeomorphically, holomophically, polynomially, formally) conjugated if there exists a germ of a homeomorphism ($C^r$ diffeomorphism, biholomorphism, polynomial biholomorphism, invertible $n$-ple of formal power series) $h: (\mathbb{C}^n,0)\rightarrow(\mathbb{C}^n,0)$, such that $h\circ f=g\circ h$. A typical problem in local dynamical systems is to find a classification of germs under the equivalence relation of conjugacy, in order to reduce to a (possibly small) class of effectively interesting cases. See ~\cite{abate} for a survey. The main source of informations in this direction is the spectrum of the differential at 0. In dimension 1, a classical example is the 

\begin{theorem}[G. K\oe nigs, \textmd{~\cite{koenigs}}]\label{koenigs}
If $f:(\mathbb{C},0)\rightarrow(\mathbb{C},0)$ is a holomorphic germ such that $df_0(z)=\lambda z$, with $|\lambda|\neq0,1$, then it is locally holomorphically conjugated with $df_0$.
\end{theorem}

If $\lambda$ is a root of the unity, the topological classification follows from the

\begin{theorem}[C. Camacho, \textmd{~\cite{camacho}}]\label{camacho} Let us assume that $f:(\mathbb{C},0)\rightarrow(\mathbb{C},0)$ is a holomorphic germ such that $df_0(z)=\lambda z$, with $\lambda^q=1$, $\lambda^p\neq1$ for $0<p<q$: then, either $f^q(z)\equiv z$, or $f$ is locally topologically conjugated to $g(z):=\lambda z+z^{kq+1}$, for some $k\geq1$. 
\end{theorem}

We remark that in the above theorem, if $f^q(z)\equiv z$, then $f$ is locally holomorphically conjugated to $g(z):=\lambda z$; the conjugacy is given by the holomorphic germ $h(z):=\sum_{i=0}^{q-1}\lambda^{-i}f^i(z)$. In several dimensions, the problem of classification grows in complexity and, with some exceptions, it is wide open. A well known result is the following partial generalization of Theorem \ref{koenigs}:

\begin{theorem}[D.M. Grobman, \textmd{~\cite{grobman}, ~\cite{grobman2}}; P. Hartman,  \textmd{~\cite{hartman}}]\label{grobman hartman}
If $f:(\mathbb{C}^n,0)\rightarrow(\mathbb{C}^n,0)$ is a holomorphic germ such that the eigenvalues of $df_0$ have modulus different from 0,1, then it is locally topologically conjugated with $df_0$.
\end{theorem}

The germ above is \emph{hyperbolic}. We shall consider the problem of topological classification of a \emph{semi-hyperbolic} germ $f:(\mathbb{C}^n,0)\rightarrow(\mathbb{C}^n,0)$:

\begin{definition}
If $f:(\mathbb{C}^n,0)\rightarrow(\mathbb{C}^n,0)$ is a holomorphic germ fixing the origin, indicating with $\lambda_1,\ldots,\lambda_n$ the eigenvalues of $df_0$, we shall say that $f$ is \emph{semi-hyperbolic} if 
\begin{enumerate}
\item[-] there is $q\geq1$ such that $\lambda_1^q=1$, $\lambda_1^p\neq1$ for $0<p<q$;
\item[-] $|\lambda_i|\neq 0,1$ for $i=2,\ldots,n$. 
\end{enumerate}
We shall say that we are in \emph{quasi-absence of resonances} if $\lambda_2^{r_2}\cdot\ldots\cdot\lambda_n^{r_n}\neq1$, for all multi-indices $(r_2,\ldots,r_n)$ such that $r_i\geq0$ for $i=2,\ldots,n$, $\sum_{i=2}^n r_i\geq1$. 
\end{definition}

We shall prove the following 

\begin{theorem}\label{normal form first coordinate}
Let us assume that $f:=(f_1,\ldots,f_n):(\mathbb{C}^n,0)\rightarrow(\mathbb{C}^n,0)$ is a holomorphic germ, and that $\lambda_1,\ldots,\lambda_n$ are the eigenvalues of $df_0$. If $f$ is semi-hyperbolic, and if there is quasi-absence of resonances, then up to a holomorphic conjugacy one of the following is satisfied:
\begin{enumerate}
\item[$(i)$] $f_1(z)=f_1(z_1,\ldots,z_n)=\lambda_1 z_1$;
\item[$(ii)$] $f_1(z)=f_1(z_1,\ldots,z_n)=\lambda_1 z_1+a_k z_1^{kq+1}+o(|z|^{kq+1})$, with $a_k\neq0$ and $k\geq1$.
\end{enumerate}
Moreover, in case $(i)$, $f$ is locally topologically conjugated at the origin with $g:=df_0$; and, in case $(ii)$, it is locally topologically conjugated at the origin with $g(z)$ $:=$ $(\lambda_1 z_1+z_1^{kq+1},$ $\lambda_2 z_2,$ $\ldots,$ $\lambda_n z_n)$.
\end{theorem}

This result, which appears natural in view of Theorems \ref{camacho}, \ref{grobman hartman}, in dimension $n=2$ was proved by J. C. Canille Martins ~\cite{canille martins}. 

The paper is entirely devoted to prove this theorem. Let us summarize as the proof is organized.
We shall use the theory related to the invariant manifolds theorems, see for example ~\cite{hirsch pugh shub}, ~\cite{katok}, ~\cite{ruelle}. Let us introduce a

\begin{definition}
If $M$ is a Riemanniann manifold, and if $f:M\rightarrow M$ is a diffeomorphism, we shall say that a submanifold $N\subseteq M$ is \emph{normally hyperbolic for $f$} if it is $f$-invariant, and if there are $0<\lambda<\lambda'\leq1\leq\mu'<\mu$, and a $df$-invariant splitting $T_NM=TN\oplus E^s(N)\oplus E^u(N)$, such that $\| df|_{E^s(N)}\|\leq\lambda$,  $\| df|_{TN}^{-1}\|\leq(\lambda')^{-1}$, $\| df|_{TN}\|\leq\mu'$,  $\| df|_{E^u(N)}^{-1}\|\leq\mu^{-1}$.
\end{definition}

In ~\cite{palis takens}, it is considered the problem of extending topological equivalence defined on normally invariant submanifolds to  neighborhoods of them. We shall work with the same plan in mind. We shall use the following theorem, which can be considered as a generalization of Theorem \ref{grobman hartman}, from a hyperbolic invariant point to a normally hyperbolic invariant submanifold:

\begin{theorem}[C. Pugh, M. Shub, \textmd{~\cite{pugh shub}}]\label{linearizzazione sul fibrato tangente}
If $M$ is a Riemannian manifold, $f:M\rightarrow M$ is a diffeomorphism, and $N\subseteq M$ is a normally hyperbolic submanifold for $f$, with $df$-invariant splitting $T_NM=TN\oplus E^s(N)\oplus E^u(N)$, then there  are open neighborhoods $U\subseteq M$ of $N$, $\mathcal{U}\subseteq E^s(N)\oplus E^u(N)$ of the image of the zero section, and a homeomorphism $h:U\rightarrow\mathcal{U}$, such that $h\circ f(z)=df\circ h(z)$, for $z\in U\cap f^{-1}(U)$. 
\end{theorem}

By the center manifold theorem (~\cite{ruelle}, p. 32), for $0\leq r<\infty$ there exists a germ of a $C^r$,  $f$-invariant, submanifold $S$ through $0\in\mathbb{C}^n$, tangent at 0 to the eigenspace of $\lambda_1$. In general, $S$ is not unique, and it cannot in general be choosen $C^\infty$. Using the quasi-absence of resonances, and proceeding as in the proof of Theorem \ref{camacho}, in Sections \ref{sec:1}, \ref{sec:2} we shall see that $f|_S$ behaves either as $z_1\mapsto\lambda_1 z_1$, or as $z_1\mapsto \lambda_1 z_1+z_1^{kq+1}$, for some $k\geq1$. Assuming, up to a $C^r$ conjugacy, that $S\subseteq\mathbb{C}_{z_1}:=\{z_2=\ldots=z_n=0\}\subseteq\mathbb{C}^n$, we shall then consider a suitable diffeomorphism $\tilde{f}:\mathbb{C}^n\rightarrow\mathbb{C}^n$ which extends a representative of $f$, such that $\mathbb{C}_{z_1}$ is normally hyperbolic for $\tilde{f}$. The linear subspace $\mathbb{C}_{z_1}$ will be normally hyperbolic also for a suitable diffeomorphism $\tilde{g}:\mathbb{C}^n\rightarrow\mathbb{C}^n$ extending a representative of the germ $g:=df_0$,  respectively $g(z_1,\ldots,z_n):=(\lambda_1 z_1+z_1^{kq+1},\lambda_2 z_2,\ldots,\lambda_n z_n)$,  with $\tilde{f}|_{\mathbb{C}_{z_1}}$ and $\tilde{g}|_{\mathbb{C}_{z_1}}$ topologically conjugated. Working in two steps, separately on the contracting and on the expanding component of the splitting, we shall see that $d\tilde{f}|_{E^s_{\tilde{f}}(\mathbb{C}_{z_1})\oplus E^u_{\tilde{f}}(\mathbb{C}_{z_1})}$ and $d\tilde{g}|_{E^s_{\tilde{g}}(\mathbb{C}_{z_1})\oplus E^u_{\tilde{g}}(\mathbb{C}_{z_1})}$ are topologically conjugated. Technicalities for the construction of such conjugacy are in Section \ref{sec:3};
essentially we shall generalize on fibers the fact that invertible linear contractions preserving the orientation of $\mathbb{C}^n$ are topologically conjugated. The assertion will then follow from Theorem \ref{linearizzazione sul fibrato tangente}.

I would like to thank M. Abate for suggesting to me the problem, and for his precious support, essential for the realization of this paper; and S. Marmi for a couple of stimulating conversations.


\section{Normal Form in the First Coordinate}
\label{sec:1}

In this section, up to a formal conjugacy, we shall reduce to a simple form the first coordinate of a holomorphic, semi-hyperbolic germ, with quasi-absence of resonances. We shall prove the following

\begin{proposition}\label{formal conjugacy}
Let us assume that $f:=(f_1,\ldots,f_n):(\mathbb{C}^n,0)\rightarrow(\mathbb{C}^n,0)$ is a holomorphic, semi-hyperbolic germ, with quasi-absence of resonances, with $\lambda_1,\ldots,\lambda_n$ eigenvalues for $df_0$. Let us assume also that $f(z)=Az+\phi(z)$, where $A$ is in Jordan form and $diag(A)=(\lambda_1,\ldots,\lambda_n)$,
with $\phi(z)=o(|z|)$. Then, there exists a formal power series $\xi_1(z)$ $=\sum_{|P|\geq2}$ $\xi_{1,P}z^P$ $\in\mathbb{C}[[z]]$ such that, if $h(z)=z+\xi(z)$ with $\xi(z):=(\xi_1(z),0,\ldots,0)$, and if $h\circ f\circ h^{-1}=(\tilde{f}_1,\ldots,\tilde{f}_n)$, then
\begin{equation}\label{forma normale prima coordinata}
\tilde{f}_1(z)=\lambda_1 z_1+\sum_{k=1}^\infty a_k z_1^{kq+1}.
\end{equation}
Moreover, if $q=1$, it is possible to assume that $\xi_1(z_1,0,\ldots,0)\equiv0$.
\end{proposition}

\emph{Proof.} 
To simplify the notations, during the proof we shall continue to use the same letter $f$ also for its conjugated maps.  We shall prove the lemma constructing a countable sequence of polynomial conjugations, and considering its limit. 
Let us assume we have reduced to the case in which 
\begin{equation}\label{normal form passo k}
\phi_1(z)=\sum_{|P|=2}^k\psi_{1,P}z^P+\sum_{|P|\geq k+1}\phi_{1,P}z^P,
\end{equation}
with $\psi_{1,P}=0$ for $P\neq(hq+1,0,\ldots,0)$, for $h\geq 1$. We are going to prove the existence of a map $h^{(k+1)}(z)=z+\xi(z)$, with $\xi(z)=(\xi_1(z),0,\ldots,0)$ and $\xi_1(z)=\sum_{|P|=k+1}\xi_{1,P}z^P$, such that $h^{(k+1)}\circ f\circ (h^{(k+1)})^{-1}(z)=Az+\tilde{\phi}(z)$, with 
\begin{equation}\label{normal form passo k+1}
\tilde{\phi}_1(z)=\sum_{|P|=2}^k\psi_{1,P}z^P+\tilde{\psi}_{1,(k+1,0,\ldots,0)}z_1^{k+1}+\sum_{|P|\geq k+2}\tilde{\phi}_{1,P}z^P;
\end{equation}
moreover, if $k$ is an integer multiple of $q$, then $\xi^{(k+1)}_{1,(k+1,0,\ldots,0)}=0$; otherwise, $\tilde{\psi}_{1,(k+1,0,\ldots,0)}=0$. Writing $(h^{(k+1)})^{-1}(z)=z+\hat{\xi}(z)$, with $\hat{\xi}(z) $ $=$ $(\hat{\xi}_1(z),0,$ $\ldots,0)$ and $\hat{\xi}_1(z)=-\xi_1(z)+o(| z|^{k+1})$, we have
\begin{eqnarray*}
 & & h^{(k+1)}\circ f\circ (h^{(k+1)})^{-1}(z) \\
 & = & h^{(k+1)}\circ f(z+\hat{\xi}(z))\\
 & = & h^{(k+1)}(Az+A\hat{\xi}(z)+\phi(z+\hat{\xi}(z)))\\
 & = & Az+A\hat{\xi}(z)+\phi(z+\hat{\xi}(z))+\xi(Az+A\hat{\xi}(z)+\phi(z+\hat{\xi}(z))).
\end{eqnarray*}
If $h^{(k+1)}\circ f\circ (h^{(k+1)})^{-1}=(\tilde{f}_1,\ldots,\tilde{f}_n)$, then
\[
\tilde{f}_1(z)=\lambda_1 z_1-\lambda_1\xi_1(z)+\sum_{|P|=2}^k\psi_{1,P}z^P+\sum_{|P|= k+1}\phi_{1,P}z^P+\xi_1(Az)+o(| z|^{k+1}).
\]
We want that 
\begin{equation}\label{condizione per la forma normale k+1}
-\lambda_1\xi_1(z)+\sum_{|P|= k+1}\phi_{1,P}z^P+\xi_1(Az)=\tilde{\psi}_{1,(k+1,0,\ldots,0)}z_1^{k+1}.
\end{equation}
It is useful to introduce the lexicographic order over $\{P=(p_1,\ldots,p_n)\,\,\big|\,\,|P|=k+1\}$:  we shall write 
$P\prec P'$ if $p_j<p_j'$, with $j:=\inf\{i=1,\ldots,n\,\,|\,\,p_i\neq p_i'\}$. Let us note that
\begin{eqnarray*}
\xi_1(Az) & = & \sum_{|P|=k+1}\xi_{1,P}(\lambda_1 z_1)^{p_1}(\lambda_2 z_2+\rho_2 z_3)^{p_2} \\
 & & \ldots(\lambda_{n-1} z_{n-1}+\rho_{n-1} z_n)^{p_{n-1}}(\lambda_n z_n)^{p_n},
\end{eqnarray*}
with $\rho_i\in\{0,1\}$ for $i=2,\ldots,n-1$. Since 
\begin{eqnarray*}
 & & (\lambda_1 z_1)^{p_1}(\lambda_2 z_2+\rho_2 z_3)^{p_2}\ldots(\lambda_{n-1} z_{n-1}+\rho_{n-1} z_n)^{p_{n-1}}(\lambda_n z_n)^{p_n}\\
 & = & \Lambda^P z^P+\sum_{|P'|=k+1, P'\prec P}C_{PP'}z^{P'}
\end{eqnarray*}
for suitable constants $C_{P,P'}$ and with $\Lambda=(\lambda_1,\ldots,\lambda_n)$, it follows that $\xi_1(Az)=\sum_{|P|=k+1}$ $\xi_{1,P}'z^P$, with
\[
\xi_{1,P}'=\Lambda^P\xi_{1,P}+\sum_{|P'|=k+1, P\prec P'}\xi_{1,P'}C_{P'P}.
\]
In particular, $\xi_{1,(k+1,0,\ldots,0)}'=\lambda_1^{k+1}\xi_{1,(k+1,0,\ldots,0)}$. Putting
\[
\xi_{1,(k+1,0,\ldots,0)}:=\left\{
\begin{array}{ll}
(\lambda_1-\lambda_1^{k+1})^{-1}\phi_{1,(k+1,0,\ldots,0)} & \mbox{if $\frac{k}{q}\notin\mathbb{N}$},\\
0 & \mbox{otherwise};
\end{array}
\right.
\]
\[
\tilde{\psi}_{1,(k+1,0,\ldots,0)}:=\left\{
\begin{array}{ll}
0 & \mbox{if $\frac{k}{q}\notin\mathbb{N}$},\\
\phi_{1,(k+1,0,\ldots,0)} & \mbox{otherwise},
\end{array}
\right.
\]
it is possible to define recursively $\xi_{1,P}$, for $|P|\prec (k+1,0,\ldots,0)$, so that (\ref{condizione per la forma normale k+1}) is satisfied, since $\lambda_1-\Lambda^P\neq 0$ for $|P|\neq (k+1,0,\ldots,0)$, $|P|=k+1$. We have hence proved the existence of a map $h^{(k+1)}$ such that, if the coefficients of the Taylor expansion of $f$ satisfy (\ref{normal form passo k}), then those of $h^{(k+1)}\circ f\circ (h^{(k+1)})^{-1}$ satisfy (\ref{normal form passo k+1}), for $k\geq1$; moreover, $h^{(k+1)}(z_1,0,\ldots,0)=(z_1,0,\ldots,0)$ if $q=1$. To conclude the proof, it is sufficient to define $h:=\lim_{k\rightarrow\infty}h^{(k)}\circ\ldots\circ h^{(2)}$, where the limit must be considered in $(\mathbb{C}[[z]])^n$. $\diamondsuit$


\section{A Topological Conjugacy in Dimension 1}
\label{sec:2}

Working in the same spirit of Theorem \ref{camacho}, in this section we shall construct a topological conjugacy for particular $C^r$ germs at $0\in\mathbb{C}$. Preliminarly, it is useful to fix some notations. If $\mathbb{C}^{\ast}:=\mathbb{C}\setminus\{0\}$, and $U_l(i):=\{w\in\mathbb{C}^\ast\,\,|\,\,-l\pi<\arg(w)<2\pi-l\pi\}$ for $l=0,1$, $i\in\mathbb{N}$, let us define
\[
\mathbb{C}^\ast[kq]:=\Big(\bigcup_{l=0}^1\bigcup_{i=0}^{kq-1}U_l(i)\Big)\bigg/ \sim,
\]
where the union is disjoint, and $\sim$ is the equivalence relation such that $w\in U_0(i)\,\sim\, w'\in U_1(j)$ iff $w=w'$ and either $0<\arg(w)<\pi$ with $i=j$, or $\pi<\arg(w)<2\pi$ with $j=(i-1) \,\,mod(kq)$. Defining the charts $\psi_l(i):U_l(i)\hookrightarrow \mathbb{C}^\ast[kq]$, it gets the structure of a 1-dimensional complex manifold. There is a biholomorphism $\gamma:\mathbb{C}^\ast\rightarrow\mathbb{C}^\ast[kq]$, defined as
\[
\gamma:z\in S_l(i)\mapsto\psi_l(i)\Big(\frac{1}{z^{kq}}\Big),
\]
where $S_l(i):=\{z\in\mathbb{C}^\ast\,\,|\,\,\frac{2\pi i+l\pi}{kq}<\arg(z)<\frac{2\pi (i+1)+l\pi}{kq}\}$, for $l=0,1$, $i=0,\ldots,kq-1$. For $R>0$, define also 
\[
\mathbb{C}_R^\ast[kq]:=\{\psi_l(i)(w)\,\,|\,\,\mbox{$w\in U_l(i)$, $|w|>R$,  $l=0,1$, $i=0,\ldots,kq-1$}\}.
\]
Let us introduce the following lemmata.

\begin{lemma}\label{traslazione all'infinito}
If $\lambda:=e^{2\pi i\frac{h}{q}}$ with $0<h\leq q$, $(h,q)=1$, $k\geq1$, and $R>kq$,  let us consider the holomorphic function $\Phi:\mathbb{C}^\ast_R[kq]\rightarrow\mathbb{C}^\ast[kq]$, such that 
\[
\Phi:\psi_l(i)(w)\mapsto\psi_l(i')(w-kq),
\]
where $i':=(i+hk)\,\,mod(kq)$. Then, if $\rho:=R^{-kq}$, the holomorphic function $\phi(z):=\gamma^{-1}\circ\Phi\circ\gamma(z)$ is well defined for $|z|<\rho$, with Taylor expansion at the origin $\phi(z)=\lambda z+\lambda z^{kq+1}+o(|z|^{kq+1})$.
\end{lemma}

\emph{Proof.} Since $\phi$ cannot be linear, we can assume that its Taylor expansion at 0 is $\phi(z)=az+bz^{m+1}+o(|z|^{m+1})$, with $a,b\in\mathbb{C}$, $b\neq0$, $m\geq1$. If $z\in S_l(i)$, $|z|<\rho$, let us consider $w:=z^{-kq}\in U_l(i)$. We have that
\begin{eqnarray*}
\psi_l(i')(w-kq) & = & \Phi(\big(\psi_l(i)(w)\big)\\
 & = & \gamma\circ\phi\circ\gamma^{-1}(\big(\psi_l(i)(w)\big)\\
 & = & \gamma\big(az+bz^{m+1}+o(|z|^{m+1})\big)\\
 & = & \psi_{l'}(j)\big((az+bz^{m+1}+o(|z|^{m+1}))^{-kq}\big)\\
 & = & \psi_{l'}(j)\big(z^{-kq}(a^{-kq}-kqba^{-(kq+1)}z^m+o(|z|^m))\big),
\end{eqnarray*}
for suitable $l'$, $j$. Necessarily, $w-kq=w(a^{-kq}-kqba^{-(kq+1)}z^m+o(|z|^m))$; hence,
$m=kq$, $a^{-kq}=1$, and $b=a$. In particular, $a=e^{2\pi i\frac{h'}{kq}}$ for some $0<h'\leq kq$. Note that $az\in S_l(i'')$, with $i''=(i+h')mod(kq)$. Since the arbitrary of $z$, we can conclude that $h'=hk$ and $a=\lambda$. $\diamondsuit$

\begin{lemma}\label{bump function}
Given $\eta>0$, there exists a $C^\infty$ function $\rho_\eta:\mathbb{C}^n\rightarrow[0,1]$, such that $\rho_\eta(z)=\rho_\eta(|z|)$, $\rho_\eta(z)=1$ for $|z|<\frac{\eta}{2}$, $supp(\rho_\eta)\subseteq\{|z|<\eta\}$, $\|d(\rho_\eta)_z\|\leq\frac{4}{\eta}$ for $z\in\mathbb{C}^n$.
\end{lemma}

We shall prove the

\begin{proposition}\label{risultato similar camacho}
In the notations introduced in Lemma \ref{traslazione all'infinito}, let us assume that $f$ is a $C^r$ function defined in a neighborhood of $0\in\mathbb{C}$, with Taylor expansion at the origin $f(z)=\lambda z+\lambda z^{kq+1}+o(|z|^{kq+1})$, where $r>kq+1$. Then, there exists $0<\eta_0<\rho$ such that, for $0<\eta<\eta_0$, the function $\tilde{f}(z):=\rho_\eta(z)f(z)+(1-\rho_\eta(z))\phi(z)$ is well defined for $|z|<\rho$, and there is a homeomorphism $\Gamma_f:\mathbb{C}\rightarrow\mathbb{C}$ such that $\Gamma_f(z)=z$ for $|z|\geq\eta$, and  $\Gamma_f\circ\phi(z)=\tilde{f}\circ\Gamma_f(z)$ for $|z|<\rho$. Here, $\rho_\eta$ is the function introduced in Lemma \ref{bump function}, with $n=1$. 
\end{proposition}

\emph{Proof.} In order to conjugate $\tilde{f}$ to $\phi$, it suffices to conjugate $\tilde{F}$ to $\Phi$, where $\tilde{F}:=\gamma\circ\tilde{f}\circ\gamma^{-1}|_{\mathbb{C}^\ast_R[kq]}$. To this aim, first of all let us  observe that $\tilde{F}$ is ``uniformly close'' to $\Phi$. To make this claim more precise,  it is useful to introduce the regions
\[
\hat{U}_l(i):=\{w\in\mathbb{C}^\ast\,\,|\,\,-l\pi+c<\arg(w)<2\pi-l\pi-c\}, 
\]
\[
\hat{S}_l(i):=\Big\{z\in\mathbb{C}^\ast\,\,|\,\,\frac{2\pi i+l\pi+c}{kq}<\arg(z)<\frac{2\pi (i+1)+l\pi-c}{kq}\Big\},
\]
where $0<c<\frac{\pi}{2}$ is fixed, and $l=0,1$, $i=0,\ldots,kq-1$. Up to choosing $\eta_0$ sufficiently small, for $w=z^{-kq}\in\hat{U}_l(i)$, with $|w|>R$ and $z\in\hat{S}_l(i)$, we have that $\psi_l(i')^{-1}\circ\tilde{F}\circ\psi_l(i)(w)$ is well defined, since $\tilde{f}\big(\hat{S}_l(i)\big)\subseteq S_l(i')$. Moreover, 
\begin{eqnarray*}
 & & |\psi_l(i')^{-1}\circ\tilde{F}\circ\psi_l(i)(w)-\psi_l(i')^{-1}\circ\Phi\circ\psi_l(i)(w)| \\
 & = & 
|\psi_l(i')^{-1}\circ\gamma\circ\tilde{f}(z)-\psi_l(i')^{-1}\circ\gamma\circ\phi(z)|\\
 & = & \bigg| \frac{1}{\big(\rho_\eta(z)f(z)+(1-\rho_\eta(z))\phi(z)\big)^{kq}}-\frac{1}{\big(\phi(z)\big)^{kq}}\bigg|\\
 & = & \bigg| \frac{\big(\phi(z)\big)^{kq}-\big(\phi(z)+\rho_\eta(z)(f(z)-\phi(z))\big)^{kq}}{\big(\phi(z)\big)^{kq}\big(\phi(z)+\rho_\eta(z)(f(z)-\phi(z))\big)^{kq}}\bigg|\\
 & = & o(|z|),
\end{eqnarray*}
since $|f(z)-\phi(z)|=o(|z|^{kq+1})$. Hence, if $0<\eta<\eta_0$ and $\hat{R}:=\eta^{-kq}$, we can assume that $\tilde{F}\equiv\Phi$ on $\mathbb{C}^\ast_R[kq]\setminus\mathbb{C}^\ast_{\hat{R}}[kq]$, and that $\tilde{F}$ is uniformly close to $\Phi$ on $\mathbb{C}^\ast_{\hat{R}}[kq]$, in the sense above stated. Let us define the region 
\[
\Omega_0:=\{z\in\mathbb{C}\,\,\big|\,\,|z|\leq\hat{R}\}\cup\{z=x+iy\in\mathbb{C}\,\,\big|\,\,0\leq x\leq kq\}\subseteq\mathbb{C},
\]
obtained by the union of a disk and a vertical strip in the complex plane; and let
\[
\Omega:=\bigcup_{l=0}^1\bigcup_{i=0}^{kq-1}\psi_l(i)\big(U_l(i)\cap\Omega_0\big)\subseteq\mathbb{C}^\ast[kq].
\]
We want to define a homeomorphism with the image $\Gamma_0:\Omega\rightarrow \mathbb{C}^\ast[kq]$, such that $\Gamma_0\circ\Phi=\tilde{F}\circ\Gamma_0$, when the compositions are well defined. To this aim, for $\zeta=\psi_l(i')(w)$ with $w=x+iy\in U_l(i')\cap\Omega_0$, put $\Gamma_0(\zeta):=\zeta$ if $|w|\leq\hat{R}$; put $\Gamma_0(\zeta):=\zeta$ if $x=kq$; put $\Gamma_0(\zeta):=\tilde{F}\circ\psi_l(i)(w+kq)$ if $x=0$, $|y|>\hat{R}$. Extend then $\Gamma_0$ to the whole $\Omega$, so that it is a homeomorphism with the image. This is possible, since $\tilde{F}$ is close to $\Phi$. The next step is to define a homeomorphism $\Gamma:\mathbb{C}^\ast[kq]\rightarrow\mathbb{C}^\ast[kq]$, as $\Gamma(\zeta):=\Gamma_0(\zeta)$ if $\zeta\in\Omega$, and $\Gamma(\zeta):=\tilde{F}^{-m}\circ\Gamma_0\circ\Phi^m(\zeta)$ if $\zeta\notin\Omega$, where $m\in\mathbb{Z}$ is the  smallest in modulus value such that $\Phi^m(\zeta)\in\Omega$. By construction, it is such that $\Gamma\circ\Phi(\zeta)=\tilde{F}\circ\Gamma(\zeta)$, for $\zeta\in\mathbb{C}^\ast_R[kq]$; and $\Gamma(\zeta)=\zeta$, for $\zeta\in\mathbb{C}^\ast[kq]\setminus\mathbb{C}^\ast_{\hat{R}}[kq]$. It suffices to define $\Gamma_f:=\gamma^{-1}\circ\Gamma\circ\gamma$ to conclude the proof. $\diamondsuit$


\section{Normally Hyperbolic Submanifolds}
\label{sec:3}

We shall now consider some results related to the theory of invariant normally hyperbolic submanifolds, useful to prove Theorem \ref{normal form first coordinate}. Let us introduce a notation. If $f:U\rightarrow\mathbb{C}^n$ is a $C^1$ function defined on  $U\subseteq\mathbb{C}^n$, we shall indicate with $\|f\|_{1,U}$ the supremum on $U$ of the norm of $f$ and of its differential. 

\begin{proposition}\label{costruzione della varieta normalmente iperbolica}
Let us assume that $L:\mathbb{C}^n\rightarrow\mathbb{C}^n$ is linear invertible, and that $\mathbb{C}^n=\mathbb{C}^h\oplus\mathbb{C}^k\oplus\mathbb{C}^l$ is a $L$-invariant splitting, such that $\| L|_{\mathbb{C}^k}\|<\lambda$, $\| L|_{\mathbb{C}^h}^{-1}\|<\lambda^{-1}$, $\| L|_{\mathbb{C}^h}\|<\mu$, $\| L|_{\mathbb{C}^l}^{-1}\|<\mu^{-1}$, for some $0<\lambda<1<\mu$. Given $\gamma>0$, there exists $0<\epsilon=\epsilon(\gamma, L)<\gamma$ such that, if $f:\mathbb{C}^n\rightarrow\mathbb{C}^n$ is a diffeomorphism, \mbox{$\| f-L\|_{1,\mathbb{C}^n}<\epsilon$}, and $\mathbb{C}^h\cong\mathbb{C}^h\oplus\{0\}\oplus\{0\}\subseteq\mathbb{C}^n$ is $f$-invariant, then $\mathbb{C}^h\subseteq\mathbb{C}^n$ is normally hyperbolic for $f$, with splitting $T_{\mathbb{C}^h}\mathbb{C}^n=T\mathbb{C}^h\oplus E^s(\mathbb{C}^h)\oplus E^u(\mathbb{C}^h)$. Moreover, identifying canonically  $T_{\mathbb{C}^h}\mathbb{C}^n\cong\mathbb{C}^h\times(\mathbb{C}^h\oplus\mathbb{C}^k\oplus\mathbb{C}^l)$, the fibers of the splitting satisfy the following relations:
\begin{eqnarray*}
E^s_x(\mathbb{C}^h)\subseteq\big\{(v_1,v_2,v_3)\in \mathbb{C}^h\oplus\mathbb{C}^k\oplus\mathbb{C}^l \,\,|\,\,|(v_1,0,v_3)|\leq\gamma| v_2|\big\},\\
E^u_x(\mathbb{C}^h)\subseteq\big\{(v_1,v_2,v_3)\in \mathbb{C}^h\oplus\mathbb{C}^k\oplus\mathbb{C}^l \,\,|\,\,|(v_1,v_2,0)|\leq\gamma| v_3|\big\},
\end{eqnarray*}
for $x\in\mathbb{C}^h$.
\end{proposition}

In order to prove it, let us introduce a definition. 
If $\mathbb{C}^n=\mathbb{C}^\alpha\oplus\mathbb{C}^\beta$, and $\gamma_0>0$,  the \emph{standard horizontal $\gamma_0$-cone at 0} is 
\[
H^{\gamma_0}=\{(u,v)\in\mathbb{C}^\alpha\oplus\mathbb{C}^\beta \,\,\big|\,\,| v|\leq\gamma_0| u|\},
\]
while the \emph{standard verical $\gamma_0$-cone at 0} is 
\[
V^{\gamma_0}=\{(u,v)\in\mathbb{C}^\alpha\oplus\mathbb{C}^\beta \,\,\big|\,\,| u|\leq\gamma_0| v|\}.
\]
Proposition \ref{costruzione della varieta normalmente iperbolica} follows as a consequence of the next, well known, lemma (a proof can be found in ~\cite{katok}, p. 248).

\begin{lemma}\label{famiglie di coni}
Assume that $\{L_m:\mathbb{C}^\alpha\oplus\mathbb{C}^\beta\rightarrow\mathbb{C}^\alpha\oplus\mathbb{C}^\beta\}_{m\in \mathbb{Z}}$ is a family of invertible linear maps, such that 
\begin{enumerate}
\item[$(i)$] $L_m(H^{\gamma_0})\subseteq H^{\gamma_0}$,
\item[$(ii)$] $L_m^{-1}(V^{\gamma_0})\subseteq V^{\gamma_0}$;
\item[$(iii)$] $| L_m(u,v)|\geq c|(u,v)|$ for $(u,v)\in H^{\gamma_0}$,
\item[$(iv)$] $| L_m(u,v)|\leq c'|(u,v)|$ for $(u,v)\in L_m^{-1}(V^{\gamma_0})$,
\end{enumerate}
for some $\gamma_0>0$ and $0<c'<c$, independent of $m$.
Then, 
\[
E^+:=\bigcap_{m=1}^\infty L_{-1}\circ\ldots\circ L_{-m}(H^{\gamma_0})\subseteq H^{\gamma_0}
\]
is a $\alpha$-dimensional linear subspace,
and 
\[
E^-:=\bigcap_{m=1}^\infty L_1^{-1}\circ\ldots\circ L_m^{-1}(V^{\gamma_0})\subseteq V^{\gamma_0}
\]
is a $\beta$-dimensional linear subspace.
\end{lemma}

\emph{Proof of Proposition \ref{costruzione della varieta normalmente iperbolica}.} To construct the contracting component of the splitting, let us reorder the basis of $\mathbb{C}^n$, so that the splitting $\mathbb{C}^n=\mathbb{C}^{h+l}\oplus\mathbb{C}^k$ is $L$-invariant, and $\| L|_{\mathbb{C}^k}\|<\lambda$,  $\| L|_{\mathbb{C}^{h+l}}^{-1}\|<\lambda^{-1}$. Define the sequence $\hat{L}_m:=L:\mathbb{C}^{h+l}\oplus\mathbb{C}^k\rightarrow\mathbb{C}^{h+l}\oplus\mathbb{C}^k$, for $m\in\mathbb{Z}$. It verifies conditions $(i)-(iv)$ of the previous lemma, with $\alpha:=h+l$, $\beta=k$, and suitable positive constants $\hat{\gamma}_0$ and $\hat{c}'<\lambda<\hat{c}$.  For $x_0\in\mathbb{C}^h\oplus\{0\}\subseteq\mathbb{C}^{h+l}$, let us consider the sequence $\{x_m\}_{m\in\mathbb{Z}}\subseteq\mathbb{C}^h\oplus\{0\}$, such that $f^m(x_0)=x_m$. Define
\[
L_m:=df_{x_m}:T_{x_m}(\mathbb{C}^n)\cong(\mathbb{C}^{h+l}\oplus\mathbb{C}^k)\rightarrow
T_{x_{m+1}}(\mathbb{C}^n)\cong(\mathbb{C}^{h+l}\oplus\mathbb{C}^k).
\]
Take $0<\epsilon=\epsilon(\gamma,L)<\gamma$. If it is suitably small, and if $\| f-L\|_{1,\mathbb{C}^n}<\epsilon$, it follows that the family $\{L_m\}_{m\in\mathbb{Z}}$ satisfies  the hypotheses of the previous lemma, for
$\alpha=h+l$, $\beta=k$, and suitable positive constants $\gamma_0<\min\{\hat{\gamma}_0,\gamma\}$, and $\hat{c}'<c'<\lambda<c<\hat{c}$, independent on $m$. Hence,
\[
E^-_{x_0}:=\bigcap_{i=1}^\infty L_0^{-1}\circ\ldots\circ L_i^{-1}(V^{\gamma_0})\subseteq T_{x_0}(\mathbb{C}^n)
\]
is a $k$-dimensional subspace, and the stable component of the splitting we are looking for is given by
$E^s(\mathbb{C}^h):=\bigcup_{x_0\in\mathbb{C}^h}E^-_{x_0}$. Indeed, by construction it is $df$-invariant, and $| df_x(v)|\leq c'| v|$, for $x\in\mathbb{C}^h$ and $v\in T_x(\mathbb{C}^n)$.
Up to shrink the constant $\epsilon$, with a similar argument we can construct the expanding component of the splitting, and conclude the proof. $\diamondsuit$

\,\,\,\,\,\,\,

We shall use also the

\begin{proposition}\label{coniugio dei differenziali}
Let us assume that $L:\mathbb{C}^n\rightarrow\mathbb{C}^n$ is linear invertible, with $L$-invariant splitting $\mathbb{C}^n=\mathbb{C}^h\oplus\mathbb{C}^k\oplus\mathbb{C}^l$. There exists $\gamma=\gamma(L)>0$ such that, if  $f_\alpha,f_\beta:\mathbb{C}^n\rightarrow\mathbb{C}^n$ are  diffeomorphisms, and 
\begin{enumerate}
\item[$(a)$] $\| f_i-L\|_{1,\mathbb{C}^n}<\gamma$, for $i=\alpha,\beta$;
\item[$(b)$] $\mathbb{C}^h\cong\mathbb{C}^h\oplus\{0\}\oplus\{0\}\subseteq\mathbb{C}^n$ is $f_i$-invariant, normally hyperbolic for $f_i$, with $df_i$-invariant splitting $T_{\mathbb{C}^h}\mathbb{C}^n=T\mathbb{C}^h\oplus (E^i)^s(\mathbb{C}^h)\oplus (E^i)^u(\mathbb{C}^h)$, for $i=\alpha,\beta$;
\item[$(c)$] the vector sub-bundle $(E^i)^{s(u)}(\mathbb{C}^h)\subseteq T_{\mathbb{C}^h}\mathbb{C}^n\cong\mathbb{C}^h\times(\mathbb{C}^h\oplus\mathbb{C}^k\oplus\mathbb{C}^l)$ admits a trivialization $(\pi^i)^{s(u)}:\mathbb{C}^h\times\mathbb{C}^{k(l)}\stackrel{\approx}{\rightarrow}(E^i)^{s(u)}(\mathbb{C}^h)\subseteq\mathbb{C}^h\times(\mathbb{C}^h\oplus\mathbb{C}^k\oplus\mathbb{C}^l)$, $(\pi^i)^{s(u)}(x,v):=(x,(\pi^i)^{s(u)}_x(v))$,  such that
\[
| (\pi^i)^{s(u)}_x(v)-\iota_x(v)|<\gamma| v|,
\]
for $x\in\mathbb{C}^h$, $v\in\mathbb{C}^{k(l)}$, and $i=\alpha,\beta$, where $\iota_x:\mathbb{C}^{k(l)}\hookrightarrow\mathbb{C}^h\oplus\mathbb{C}^k\oplus\mathbb{C}^l$ is the inclusion;
\item[$(d)$] there is a homeomorphism $\Gamma:\mathbb{C}^h\rightarrow\mathbb{C}^h$ such that $\Gamma\circ f_\alpha|_{\mathbb{C}^h}=f_\beta\circ\Gamma$, 
\end{enumerate}
then there exists a homeomorphism 
$h:(E^\alpha)^s(\mathbb{C}^h)\oplus (E^\alpha)^u(\mathbb{C}^h)$ $\rightarrow$ $(E^\beta)^s(\mathbb{C}^h)$ $\oplus (E^\beta)^u(\mathbb{C}^h)$,
lifting $\Gamma$,  such that $h\circ df_\alpha|_{(E^\alpha)^s(\mathbb{C}^h)\oplus (E^\alpha)^u(\mathbb{C}^h)}=df_\beta\circ h$. 
\end{proposition}

In order to prove this proposition, we shall proceed in two steps: at first, lifting the homeomorphism $\Gamma$ to a vector bundles isomorphism $h^s:(E^\alpha)^s(\mathbb{C}^h)\rightarrow (E^\beta)^s(\mathbb{C}^h)$ defined only on the contracting component of the splittig, then using it to construct the general isomorphism needed. Since the argument used in the two steps is similar, we prefer to isolate it in the following

\begin{lemma}\label{contrazione caso generale}
Let us assume that $L_0:\mathbb{C}^n\rightarrow\mathbb{C}^n$ is linear invertible, with $L_0$-invariant decomposition $\mathbb{C}^n=\mathbb{C}^p\oplus\mathbb{C}^q\oplus\mathbb{C}^r$; that $c:=2(1+\| L_0\|)$, and that
\begin{equation}\label{condizione su gamma}
0<\gamma<\min\Big\{1,\frac{1}{\| L_0^{-1}\|}\cdot\frac{1}{\| L_0\|+c}\Big\}.
\end{equation}
If\, $\mathbb{C}^p\times(\mathbb{C}^p\oplus\mathbb{C}^q\oplus\mathbb{C}^r)\rightarrow\mathbb{C}^p$ is the trivial bundle, let us assume that $E^i\subseteq\mathbb{C}^p\times(\mathbb{C}^p\oplus\mathbb{C}^q\oplus\mathbb{C}^r)$ is a rank-$r$ trivial sub-bundle, endowed with a trivialization $\pi^i:\mathbb{C}^p\times\mathbb{C}^r\rightarrow E^i\subseteq\mathbb{C}^p\times(\mathbb{C}^p\oplus\mathbb{C}^q\oplus\mathbb{C}^r)$, $\pi^i(x,v):=(x,\pi^i_x(v))$, such that
\begin{equation}\label{descrizione fibra nel caso-calco}
|\pi^i_x(v)-\iota_x(v)|<\gamma| v|,
\end{equation}
for $x\in\mathbb{C}^p$, $v\in\mathbb{C}^r$, $i=\alpha,\beta$,  where $\iota_x:\mathbb{C}^r\hookrightarrow\mathbb{C}^p\oplus\mathbb{C}^q\oplus\mathbb{C}^r$ is the inclusion. Let us also assume that $F_i:\mathbb{C}^p\rightarrow\mathbb{C}^p$ is a  homeomorphism; and that  $\phi^i:E^i\rightarrow E^i$ is a vector bundles isomorphism, lifting $F_i$, uniformly contracting (expanding) the fibers, such that 
\begin{equation}\label{df vicina ad L}
| (\phi^i)_{F_i^{-1}(x)}(\pi^i_{F_i^{-1}(x)}(v))-\pi^i_x(L_0(v))|\leq c\gamma| v|,
\end{equation}
for $x\in\mathbb{C}^p$, $v\in\mathbb{C}^r$, $i=\alpha,\beta$. Then, if \,$\Gamma_0:\mathbb{C}^p\rightarrow\mathbb{C}^p$ is a homeomorphism such that $\Gamma_0\circ F_\alpha=F_\beta\circ\Gamma_0$, there exists a homeomorphism $H:E^\alpha\rightarrow E^\beta$, lifting $\Gamma_0$, such that $H\circ\phi^\alpha=\phi^\beta\circ H$.
\end{lemma}

Assuming this lemma, let us see the

\,\,\,\,\,\,

\emph{proof of Proposition \ref{coniugio dei differenziali}.} \emph{Step 1.} In order to be in the hypotheses of Lemma \ref{contrazione caso generale}, let us put $p:=h$; $q:=l$; $r:=k$; $L_0:=L$; $E^i:=(E^i)^s(\mathbb{C}^h)$; $\pi^i:=(\pi^i)^s$; $F_i:=f_i|_{\mathbb{C}^h}$; $\phi^i:=df_i|_{(E^i)^s(\mathbb{C}^h)}$; $\Gamma_0:=\Gamma$. We have to check that the relation (\ref{df vicina ad L}) is verified. To this aim, for $x\in\mathbb{C}^h$ and $v\in\mathbb{C}^k$, we have
\begin{eqnarray*}
 & & | d(f_i)_{f_i^{-1}(x)}\big((\pi^i)^s_{f_i^{-1}(x)}(v)\big)-(\pi^i)^s_x\big(L(v)\big)|\\
 &\leq & | d(f_i)_{f_i^{-1}(x)}\big((\pi^i)^s_{f_i^{-1}(x)}(v)\big)-L\big((\pi^i)^s_{f_i^{-1}(x)}(v)\big)| \\
 & & +\,\,| L\big((\pi^i)^s_{f_i^{-1}(x)}(v)\big)-L(v)| \\
 & & +\,\,|(\pi^i)^s_x\big(L(v)\big)-L(v)|\\
 & = & (A)+(B)+(C).
\end{eqnarray*}
Now,
\begin{eqnarray*}
(A) & \leq & \gamma| (\pi^i)^s_{f_i^{-1}(x)}(v) |
 \leq\gamma\Big(|(\pi^i)^s_{f_i^{-1}(x)}(v)-v|+| v|\Big)\leq \gamma(1+\gamma)| v|,
\end{eqnarray*}
for the hypotheses $(a)$, $(c)$; similarly, 
\begin{eqnarray*}
(B) & \leq &\| L\| \cdot | (\pi^i)^s_{f_i^{-1}(x)}(v)-v|\leq\gamma\| L\|\cdot| v|;
\end{eqnarray*}
and
\[
(C)\leq \gamma\| L\|\cdot| v|.
\]
It follows that
\[
| d(f_i)_{f_i^{-1}(x)}\big((\pi^i)^s_{f_i^{-1}(x)}(v)\big)-(\pi^i)^s_x\big(L(v)\big)|\leq\gamma(1+\gamma+2\| L\|)| v|\leq c\gamma | v|.
\]
We are in the hypotheses of Lemma \ref{contrazione caso generale}. There is hence a homeomorphism $h^s:(E^\alpha)^s\rightarrow(E^\beta)^s$, lifting $\Gamma$, such that $h^s\circ df_\alpha|_{(E^\alpha)^s(\mathbb{C}^h)}=df_\beta\circ h^s$.

\emph{Step 2.} Let us consider the trivial, rank-$(k+l)$, vector bundle $P_i:(E^i)^s(\mathbb{C}^h)\oplus(E^i)^u(\mathbb{C}^h)\rightarrow\mathbb{C}^h$, with trivialization $\tilde{\pi}^i:\mathbb{C}^h\times(\mathbb{C}^k\oplus\mathbb{C}^l)\stackrel{\approx}{\rightarrow}(E^i)^s(\mathbb{C}^h)\oplus(E^i)^u(\mathbb{C}^h)$, $\tilde{\pi}^i(x,v+w):=\big(x,(\pi^i)^s_x(v)+(\pi^i)^u_x(w)\big)$. Let us define the trivial, rank-$l$, vector bundle $P_i':(E^i)^s(\mathbb{C}^h)\oplus(E^i)^u(\mathbb{C}^h)\rightarrow(\mathbb{C}^h\times\mathbb{C}^k)$, with $P_i'$ defined so that a trivialization is $(\tilde{\pi}^i)':(\mathbb{C}^h\times\mathbb{C}^k)\times\mathbb{C}^l\stackrel{\approx}{\rightarrow}(E^i)^s(\mathbb{C}^h)\oplus(E^i)^u(\mathbb{C}^h)$, where $(\tilde{\pi}^i)'\big((x,v),w\big):=\big(x,(\pi^i)^s_x(v)+(\pi^i)^u_x(w)\big)$. If $(\mathbb{C}^h\times\mathbb{C}^k)\times\big((\mathbb{C}^h\times\mathbb{C}^k)\oplus\{0\}\oplus\mathbb{C}^l\big)\rightarrow(\mathbb{C}^h\times\mathbb{C}^k)$ is the trivial bundle, let us consider the trivial sub-bundle $\hat{P}_i:\hat{E}^i\subseteq(\mathbb{C}^h\times\mathbb{C}^k)\times\big((\mathbb{C}^h\times\mathbb{C}^k)\oplus\{0\}\oplus\mathbb{C}^l\big)\rightarrow (\mathbb{C}^h\times\mathbb{C}^k)$, with trivialization $\hat{\pi}^i:(\mathbb{C}^h\times\mathbb{C}^k)\times\mathbb{C}^l\rightarrow\hat{E}^i$, $\hat{\pi}^i\big((x,v),w\big):=\big((x,v),(\pi^i)^u_x(w)\big)$. Given the expressions of the trivializations, there exists a vector bundles isomorphism $\Delta^i:(E^i)^s(\mathbb{C}^h)\oplus(E^i)^u(\mathbb{C}^h)\rightarrow\hat{E}^i$, lifting $id_{\mathbb{C}^h\times\mathbb{C}^k}$. In order to apply Lemma \ref{contrazione caso generale}, let us define this time $p:=h+k$, $q:=0$, $r:=l$; $L_0:=L$; $E^i:=\hat{E}^i$; $\pi^i:=\hat{\pi}^i$; $F_i:=\big((\pi^i)^s\big)^{-1}\circ df_i\circ(\pi^i)^s$; $\phi_i:=\Delta^i\circ df_i\circ(\Delta^i)^{-1}$; $\Gamma_0:=\big((\pi^\beta)^s\big)^{-1}\circ h^s\circ(\pi^\alpha)^s$, where $h^s$ is the homeomorphism introduced in the previous step. See Fig. \ref{fig:1}.
Proceeding as above, we obtain
\[
| (\phi^i)_{F_i^{-1}(x,v)}(\pi^i_{F_i^{-1}(x,v)}(w))-\pi^i_{(x,v)}(L_0(w))|\leq \gamma(1+\gamma+2\| L_0\|)| w|\leq c\gamma| w|,
\]
for $(x,v)\in\mathbb{C}^p$, $w\in\mathbb{C}^r$. We are in the hypotheses of Lemma \ref{contrazione caso generale}. There exists hence a homeomorphism $H:\hat{E}^\alpha\rightarrow\hat{E}^\beta$, lifting $\Gamma_0$, such that $H\circ\phi^\alpha=\phi^\beta\circ H$. Defining $h:=(\Delta^\beta)^{-1}\circ H\circ\Delta^\alpha:(E^\alpha)^s(\mathbb{C}^h)\oplus(E^\alpha)^u(\mathbb{C}^h)\rightarrow (E^\beta)^s(\mathbb{C}^h)\oplus(E^\beta)^u(\mathbb{C}^h)$, we obtain the lifting of $\Gamma$ needed to conclude the proof. $\diamondsuit$

\,\,\,\,\,\,

\begin{figure}[h]
\[
\xymatrix@1{
 & *\txt{$(E^\alpha)^s(\mathbb{C}^h)$\\
		$\oplus$ \\ $(E^\alpha)^u(\mathbb{C}^h)$} \ar[rr]^{h} \ar[dl]_{df_\alpha} \ar@{-}[d]_{\Delta^\alpha} & & 
 *\txt{$(E^\beta)^s(\mathbb{C}^h)$\\
		$\oplus$ \\ $(E^\beta)^u(\mathbb{C}^h)$} \ar[dl]_{df_\beta} \ar[dd]^{\Delta^\beta} \\
*\txt{$(E^\alpha)^s(\mathbb{C}^h)$\\
		$\oplus$ \\ $(E^\alpha)^u(\mathbb{C}^h)$} \ar[rr]^{
\,\,\,\,\,\,\,\,\,\,\,\,\,\,\,\,
\,\,\,\,\,\,\,\,\,\,\,\,\,
h} \ar[dd]_{\Delta^\alpha} & \ar[d] & 
*\txt{$(E^\beta)^s(\mathbb{C}^h)$\\
		$\oplus$ \\ $(E^\beta)^u(\mathbb{C}^h)$} \ar@{-}[d]^{\,\Delta^\beta} \ar[dd]  &  \\
 & E^\alpha \ar@{-}[r]^{H} \ar[dl]_{\phi_\alpha} \ar@{-}[d]_{\hat{P}_\alpha}& \,\,\,\,\,\ar[r] & E^\beta \ar[dl]_{\phi_\beta} \ar[dd]^{\hat{P}_\beta} &  \\
E^\alpha \ar[rr]^{
\,\,\,\,\,\,\,\,\,\,\,\,\,\,\,\,
\,\,\,\,\,\,\,\,\,\,\,\,\,\,\,\,\,\,\,\,\,\,\,\,\,\,\,\,\,\,\,\,\,\,\,
H} \ar[dd]_{\hat{P}_\alpha} & \ar[d] & E^\beta \ar@{-}[d]^{\hat{P}_\beta} \ar[dd] & \\  
 & \mathbb{C}^h\oplus\mathbb{C}^k \ar@{-}[r]^{\Gamma_0} \ar[dl]_{F_\alpha} \ar@{-}[d]_{(\pi^\alpha)^s}& \,\,\,\,\,\ar[r] & 
 \mathbb{C}^h\oplus\mathbb{C}^k \ar[dl]_{F_\beta} \ar[dd]^{(\pi^\beta)^s} &  \\
\mathbb{C}^h\oplus\mathbb{C}^k \ar[rr]^{
\,\,\,\,\,\,\,\,\,\,\,\,\,\,\,\,
\,\,\,\,\,\,\,\,\,\,\,\,\,\,\,\,\,\,\,\,\,\,\,\,\,\,\,\,\,\,\,\,\,\,\,
\Gamma_0}
 \ar[dd]_{(\pi^\alpha)^s} & \ar[d] & 
\mathbb{C}^h\oplus\mathbb{C}^k \ar@{-}[d]^{(\pi^\beta)^s} \ar[dd] & \\  
 & (E^\alpha)^s(\mathbb{C}^h) \ar@{-}[r]^{h^s} \ar[dl]_{df_\alpha} \ar@{-}[d]& \,\,\,\,\,\ar[r] & 
 (E^\beta)^s(\mathbb{C}^h) \ar[dl]_{df_\beta} \ar[dd] &  \\
(E^\alpha)^s(\mathbb{C}^h) \ar[rr]^{
\,\,\,\,\,\,\,\,\,\,\,\,\,\,\,\,
\,\,\,\,\,\,\,\,\,\,\,\,\,\,\,\,\,\,\,\,\,\,\,\,\,\,\,\,\,\,\,\,\,\,\,
h^s} \ar[dd] & \ar[d] & 
(E^\beta)^s(\mathbb{C}^h) \ar[dd] & \\  
 & \mathbb{C}^h \ar@{-}[r]^{\Gamma} \ar[dl]_{f_\alpha}& \,\,\,\,\,\ar[r] & \mathbb{C}^h \ar[dl]_{f_\beta} &  \\
\mathbb{C}^h \ar[rr]^{
\,\,\,\,\,\,\,\,\,\,\,\,\,\,\,\,
\,\,\,\,\,\,\,\,\,\,\,\,\,\,\,\,\,\,\,\,\,\,\,\,\,\,\,\,\,\,\,\,\,\,\,
\Gamma}
 & & \mathbb{C}^h & }
\]
\caption{}
\label{fig:1}
\end{figure}

Lemma \ref{contrazione caso generale} is a corollary of 

\begin{lemma}\label{coniugio caso generale}
Under the hypotheses of Lemma \ref{contrazione caso generale}, let us define $E^i(1):=\{\nu\in E^i\,\,|\,\,|\nu|< 1\}$, and
$\Omega_i:=\overline{E^i(1)\setminus\phi^i(E^i(1))}$ ($\Omega_i:=\overline{E^i(1)\setminus(\phi^i)^{-1}(E^i(1))}$), for $i=\alpha,\beta$. Then, there exists a homeomorphism $H_0:\Omega_1\rightarrow\Omega_2$, lifting $\Gamma_0$, such that
\[
H_0\circ \phi^\alpha_x(\nu_x)=\phi^\beta_{\Gamma_0(x)}\circ H_0(\nu_x) \,\,\,\,\,\,\,\,\big(H_0\circ (\phi^\alpha)_x^{-1}(\nu_x)=(\phi^\beta)_{\Gamma_0(x)}^{-1}\circ H_0(\nu_x)\big),
\]
for $x\in\mathbb{C}^p$, $\nu_x\in E^i_x$, $| \nu_x|=1$.
\end{lemma}

\,\,\,\,\,\,\,

Let us see as Lemma \ref{contrazione caso generale} follows from it.

\,\,\,\,\,\,

\emph{Proof of Lemma  \ref{contrazione caso generale}.} Let us assume that $\phi^i$ is contracting, the expanding case is analogous. For $x\in\mathbb{C}^p$ and $\nu_x\in E^\alpha_x$, $| \nu_x|\neq0$, if $m\in\mathbb{Z}$ is the smallest in modulus value such that $(\phi^\alpha)^m(\nu_x)\in\Omega_1$, let us define $H(\nu_x):=(\phi^\beta)^{-m}\circ H_0\circ(\phi^\alpha)^m(\nu_x)$, where $H_0$ is the homeomorphism introduced in Lemma \ref{coniugio caso generale}. Let us define also $H(0_x):=0_{\Gamma_0(x)}$. The function $H:E^\alpha\rightarrow E^\beta$ so defined is by construction a homeomorphism, lifting $\Gamma_0$, such that $H\circ\phi^\alpha=\phi^\beta\circ H$. $\diamondsuit$

\,\,\,\,\,\,

Before the proof of Lemma \ref{coniugio caso generale}, let us see a general property of invertible linear maps.

\begin{lemma}\label{omotopia con la identita}
If $L_0':\mathbb{C}^n\rightarrow\mathbb{C}^n$ is an invertible linear map, then there is a homotopy of invertible linear maps $\mathcal{F}:\mathbb{C}^n\times[0,1]\rightarrow\mathbb{C}^n$, such that $\mathcal{F}_0=L_0'$ and $\mathcal{F}_1=a$, with $a^2=id_{\mathbb{C}^n}$. 
\end{lemma}

\emph{Proof.} We can assume that $L_0'(v)=Tv$, where $T=(\tau_{i,j})_{i,j=1,\ldots,n}$ is a triangular matrix. If $T(t):=(t^{j-i}\cdot \tau_{i,j})_{i,j=1,\ldots,n}$, with $t\in[0,1]$, up to the homotopy of invertible linear maps $\mathcal{F}_t'(v):=T(t)v$, we can reduce to the case in which $T$ is diagonal. Moreover, we can assume that $T=diag(\lambda_1,\ldots,\lambda_n)$, with $\lambda_i$ lying on the real negative semi-axis for $i=1,\ldots,s$, and $\lambda_i$ lying on the complement in $\mathbb{C}$ for $i=s+1,\ldots,n$, for some $0\leq s\leq n$. Let us define $a(v):=Av$, with $A:=diag(a_1,\ldots,a_n)$, where $a_i:=-1$ for $i=1,\ldots,s$, and $a_i:=1$ for $i=s+1,\ldots,n$. It is such that $a^2=id_{\mathbb{C}^n}$. It suffices to define $\mathcal{F}_t''(v):=(1-t)Tv+tAv$ to conclude to proof. $\diamondsuit$

\,\,\,\,\,\,

Let us remark that the claim in the previous lemma is false for $a:=id_{\mathbb{C}^n}$. We can conclude with the

\,\,\,\,\,\,\,

\emph{proof of Lemma \ref{coniugio caso generale}.} We can assume that $\phi^i$ contracts uniformly the fibers.
It suffices to prove the existence of a homeomorpism
$\chi_i:\mathbb{C}^p\times S^{2r}\times [0,1]\rightarrow\Omega_i$, with
\begin{equation}\label{funzione sul bordo}
\chi_i(x,v,0)=\Big(x,\frac{\pi^i_x(a(v))}{| \pi^i_x(a(v))|}\Big),
\,\,
\chi_i(x,v,1)=\bigg(x,\phi^i_{F_i^{-1}(x)}\bigg(\frac{\pi^i_{F_i^{-1}(x)}(v)}{| \pi^i_{F_i^{-1}(x)}(v)|}\bigg)\bigg),
\end{equation}
where $S^{2r}=\{x\in\mathbb{C}^r \,\,|\,\,| x|=1\}$, $a$ is the function introduced in the previous lemma, and $i=\alpha,\beta$. Assuming it, and defining the homeomorphism $\hat{\Gamma}_0:\mathbb{C}^p\times S^{2r}\times [0,1]\rightarrow \mathbb{C}^p\times S^{2r}\times [0,1]$, $\hat{\Gamma}_0(x,v,t):=(\Gamma_0(x),v,t)$, it suffices to put $H_0:=\chi_\beta\circ\hat{\Gamma}_0\circ\chi_\alpha^{-1}:\Omega^\alpha\rightarrow\Omega^\beta$ to conclude the proof. Indeed, it is a homeomorphism; it lifts $\Gamma_0$; and, if $x\in\mathbb{C}^p$, $v\in S^{2r}$, and
\[
\nu_x=\frac{\pi^\alpha_x(v)}{| \pi^\alpha_x(v)|}=\frac{\pi^\alpha_x(a^2(v))}{| \pi^\alpha_x(a^2(v))|}\in E^\alpha_x,
\]
it is such that
\begin{eqnarray*}
H_0\big(\phi^\alpha_x(\nu_x)\big) & = & \chi_\beta\circ\hat{\Gamma}\circ\chi_\alpha^{-1}\Big(\phi^\alpha_x\Big(\frac{\pi^\alpha_x(v)}{| \pi^\alpha_x(v)|}\Big)\Big)\\
& = & \chi_\beta\circ\hat{\Gamma}(F_\alpha(x),v,1)\\
& = & \chi_\beta(\Gamma_0\circ F_\alpha(x),v,1)\\
& = & \phi^\beta_{\Gamma_0(x)}\bigg(\frac{\pi^\beta_{\Gamma_0(x)}(v)}{| \pi^\beta_{\Gamma_0(x)}(v)|}\bigg)\\
& = & \phi^\beta_{\Gamma_0(x)}\circ\chi_\beta(\Gamma_0(x),a(v),0)\\
& = & \phi^\beta_{\Gamma_0(x)}\circ\chi_\beta\circ\hat{\Gamma}_0(x,a(v),0)\\
& = & \phi^\beta_{\Gamma_0(x)}\circ\chi_\beta\circ\hat{\Gamma}_0\circ\chi_\alpha^{-1}\circ\chi_\alpha(x,a(v),0)\\
& = & \phi^\beta_{\Gamma_0(x)}\big(H_0(\nu_x)\big).
\end{eqnarray*}
In order to prove the existence of $\chi_i$, let us consider the vector bundles isomorphisms which lift $id_{\mathbb{C}^p}$, $\hat{a}^i:E^i\rightarrow E^i$, with $\hat{a}^i_x(\pi^i_x(v)):=\pi^i_x(a(v))$, and $\hat{\phi}^i:E^i\rightarrow E^i$, with $\hat{\phi}^i_x(\pi^i_x(v)):=\phi^i_{F_i^{-1}(x)}(\pi^i_{F_i^{-1}(x)}(v))$. There exists a homotopy of vector bundles isomorphisms lifting $id_{\mathbb{C}^p}$, $\mathcal{F}^i:E^i\times[0,1]\rightarrow E^i$, such that $\mathcal{F}^i_0=\hat{a}^i$, and $\mathcal{F}^i_1=\hat{\phi}^i$. We can construct it in two steps. Defining $\hat{L}_0^i:E^i\rightarrow E^i$, $(\hat{L}_0^i)_x(\pi^i_x(v)):=\pi^i_x(L_0(v))$, by Lemma \ref{omotopia con la identita} there exists a homotopy of vector bundles isomorphisms lifting $id_{\mathbb{C}^p}$, $(\mathcal{F}^i)'$, such that $(\mathcal{F}^i)_0'=\hat{a}^i$, $(\mathcal{F}^i)_1'=\hat{L}_0^i$. There exists a similar homotopy, $(\mathcal{F}^i)'':=(1-t)\hat{L}_0^i+t\hat{\phi}^i$, such that $(\mathcal{F}^i)_0''=\hat{L}_0^i$, $(\mathcal{F}^i)_1''=\hat{\phi}^i$. It is well defined, since for $t\in[0,1]$, $x\in\mathbb{C}^p$, $v\in\mathbb{C}^r\setminus\{0\}$, and $\nu_x=\pi^i_x(v)\in E^i_x$, we have
\begin{eqnarray*}
| \big((\mathcal{F}^i)_t''\big)_x(\nu_x)| & = & | (1-t)\pi^i_x(L_0(x))+t\phi^i_{F_i^{-1}(x)}\big(\pi^i_{F_i^{-1}(x)}(v)\big)|\\
& \geq & | \pi^i_x(L_0(v))| -c\gamma| v|\\
& \geq & | L_0(v)| -|\pi^i_x(L_0(v))-L_0(v)|-c\gamma| v|\\
& \geq & (\| L_0^{-1}\|^{-1}-\gamma\| L_0\|-c\gamma)| v|>0,
\end{eqnarray*}
where in the first inequality we have used (\ref{df vicina ad L}),  in the third (\ref{descrizione fibra nel caso-calco}), and in the last (\ref{condizione su gamma}). The homotopy $\mathcal{F}^i$ can then be obtained composing $(\mathcal{F}^i)'$ and $(\mathcal{F}^i)''$. Since $\phi^i$ contracts uniformly the fibers, there exists $\delta>0$ such that $| \phi^i_x(\nu_x)|<(1-\delta)| \nu_x|$, for $x\in\mathbb{C}^p$ and $\nu_x\in E^i_x$. It follows that the homeomorphism $\chi_i:\mathbb{C}^p\times S^{2n}\times[0,1]\rightarrow \Omega_i$  needed to conclude the proof can be defined as
\begin{eqnarray*}
 & & \chi_i(x,v,t)\\
 & := & \left\{
\begin{array}{ll}
(1-2t\delta)\frac{(\mathcal{F}^i_{2t})_x\big(\pi^i_x(v)\big)}{| (\mathcal{F}^i_{2t})_x\big(\pi^i_x(v)\big)|} & \mbox{if $t\in[0,1/2]$;}\\
\phi^i_{F_i^{-1}(x)}\bigg(\frac{(2t-1)\pi^i_{F_i^{-1}(x)}(v)}{| \pi^i_{F_i^{-1}(x)}(v)|}+\frac{(2-2t)(1-\delta)\pi^i_{F_i^{-1}(x)}(v)}{| \phi^i_{F_i^{-1}(x)}\big(\pi^i_{F_i^{-1}(x)}(v)\big)|}\bigg) & \mbox{if $t\in[1/2,1]$. $\diamondsuit$ }
\end{array}
\right.
\end{eqnarray*}


\section{Proof of the Main Theorem}
\label{sec:4}

We have all the tools we need for the

\,\,\,\,\,\,

\emph{proof of Theorem \ref{normal form first coordinate}.} First of all, up to a linear conjugacy, we can assume that $f(z)=Az+\phi(z)$, where $A$ is in Jordan form and $diag(A)=(\lambda_1,\ldots,\lambda_n)$, with $\phi=o(| z|)$. Up to considering $h\circ f\circ h^{-1}$, with 
$h(z_1,\ldots,z_n):=(z_1, z_2+\epsilon_2 z_1^2, \ldots, z_n +\epsilon_n z_1^2)$ and suitable constants $\epsilon_2, \ldots,\epsilon_n\in\mathbb{C}$,
we can assume that $(f^q)_j(z_1,0,\ldots,0)\not\equiv 0$ for $j=1,\ldots,n$. By Proposition \ref{formal conjugacy}, there exists a formal power series $\xi_1(z)\in\mathbb{C}[[z]]$ such that, putting $h(z)=z+\xi(z)$ with $\xi(z):=(\xi_1(z),0,\ldots,0)$, and $\tilde{f}=(\tilde{f}_1,\ldots,\tilde{f}_n):=h\circ f\circ h^{-1}$, the power series $\tilde{f}_1$ satisfies (\ref{forma normale prima coordinata}). If the coefficients $a_k$ introduced in (\ref{forma normale prima coordinata}) are not all vanishing, cutting off $h$ to a suitable high order we can reduce to the case $(ii)$. If the $a_k$ are all vanishing, we shall distinguish two cases. Firstly, let us assume that $q=1$. By hypothesis, $h_1\circ f(z)= h_1(z)$; in particular, $\phi_1(z)+\xi_1(f(z))=\xi_1(z)$. From Lemma \ref{formal conjugacy} we know that $\xi_1(z_1,0,\ldots,0)\equiv0$, hence $\phi_1(z_1,0,\ldots,0)+\xi_1(f(z_1,0,\ldots,0))\equiv0$. Since by assumption $f_j(z_1,0,\ldots,0)\not\equiv0$ for $j=1,\ldots,n$, we can choose $z_1$ in a neighborhood of $0\in\mathbb{C}$ such that $f(z_1,0,\ldots,0)=(\zeta_1,\ldots,\zeta_n)$, with $\zeta_i\in\mathbb{C}\setminus\{0\}$, for $i=1,\ldots,n$. It follows that $\xi_1(\zeta_1,\ldots,\zeta_n)=-\phi_1(z_1,0,\ldots,0)$. By Abel's lemma, then, the power series $\xi_1(z)$ is absolutely and uniformly convergent on compact subsets of the open set
\[
\Omega_{\zeta_1,\ldots,\zeta_n}:=\{(t_1 \zeta_1,\ldots, t_n \zeta_n) \,\,\big|\,\, t_i\in\mathbb{C}, \,\,|t_i|<1\}\subseteq\mathbb{C}^n,
\]
and we are in case $(i)$. Let us then assume that $q>1$. By the previous case, up to a holomorphic conjugacy, we can assume that $(f^q)_1(z)=z_1$. Let us introduce the map
\[
h(z):=z+\frac{f(z)}{\lambda_1}+\ldots+\frac{f^{q-1}(z)}{\lambda_1^{q-1}}.
\]
It is locally invertible in a neighborhood of the origin, since $dh_0(z)=\hat{A}z$, where $\hat{A}$ is a triangular matrix whose $i^{th}$ coefficient of the principal diagonal is $\eta_i:=\sum_{j=0}^{q-1}(\frac{\lambda_i}{\lambda_1})^j\neq0$. By construction, it is such that $h_1\circ f(z)=\lambda_1 h_1(z)$. Hence, considering $h\circ f\circ h^{-1}$, we are in case $(i)$, and the first part of the theorem is proved.

Next, we want to prove the claim about the topological classification. By the center manifold theorem, there exists a $C^r$, local submanifold $0\in S\subseteq\mathbb{C}^n$, locally $f$-invariant at 0, tangent to $\mathbb{C}_{z_1}:=\{(z_1,0,\ldots,0)\,\,|\,\,z_1\in\mathbb{C}\}\subseteq\mathbb{C}^n$ at 0, where $r>1$ in case $(i)$, $r>kq+1$ in case $(ii)$. We can assume that it is parametrized as $S:=\{(z_1,u(z_1))\,\,|\,\,z_1\in U\}$, where $U\subseteq\mathbb{C}$ is a neighborhood of the origin, and $u:U\rightarrow\mathbb{C}^{n-1}$ is a $C^r$ function, such that $u(0)=0$, $du_0\equiv 0$. Up to a $C^r$ conjugacy, by the $C^r$ diffeomorphism $h(z_1,\ldots,z_n):=(z_1,z_2+u_2(z_1),\ldots,z_n+u_n(z_1))$, we can reduce to the case $u\equiv0$ and $S\subseteq\mathbb{C}_{z_1}$. Let us consider separately the cases $(i)$, $(ii)$. In case $(i)$, if $U\subseteq\mathbb{C}_{z_1}$ is a small neighborhood of the origin, we have that $f|_{U}(z_1)=\lambda_1 z_1$. 
If $L:=df_0$, let us consider the constants $\gamma=\gamma(L)$, $\epsilon=\epsilon(L,\gamma)$ introduced in Propositions \ref{costruzione della varieta normalmente iperbolica}, \ref{coniugio dei differenziali}. Up to choosing $\eta>0$ suitably small, the function $f_\alpha(z):=\rho_\eta(z)f(z)+(1-\rho_\eta(z))L(z):\mathbb{C}^n\rightarrow\mathbb{C}^n$ is a $C^r$ diffeomorphism, such that $\| f_\alpha-L\|_{1,\mathbb{C}^n}<\epsilon$, where $\rho_\eta$ is the function introduced in Lemma \ref{bump function}. Moreover, $\mathbb{C}_{z_1}\subseteq\mathbb{C}^n$ is $f_\alpha$-invariant, hence by Proposition  \ref{costruzione della varieta normalmente iperbolica} it is normally hyperbolic for $f_\alpha$, with splitting $T_{\mathbb{C}_{z_1}}\mathbb{C}^n=T\mathbb{C}_{z_1}\oplus(E^\alpha)^s(\mathbb{C}_{z_1})\oplus(E^\alpha)^u(\mathbb{C}_{z_1})$. Defining $f_\beta:=L$, and $\Gamma:=id_{\mathbb{C}_{z_1}}$, we are in the hypotheses of Proposition \ref{coniugio dei differenziali} (condition $(c)$ for $(E^\alpha)^{s(u)}(\mathbb{C}_{z_1})$ follows from Proposition \ref{costruzione della varieta normalmente iperbolica}). There exists hence a homeomorphism $h:$ $(E^\alpha)^s(\mathbb{C}_{z_1})$ $\oplus $ $(E^\alpha)^u(\mathbb{C}_{z_1})$ $\rightarrow$ $(E^\beta)^s(\mathbb{C}_{z_1})$ $\oplus$ $(E^\beta)^u(\mathbb{C}_{z_1})$, which conjugates $df_\alpha|_{(E^\alpha)^s(\mathbb{C}_{z_1})\oplus(E^\alpha)^u(\mathbb{C}_{z_1})}$ to $df_\beta|_{(E^\beta)^s(\mathbb{C}_{z_1})\oplus(E^\beta)^u(\mathbb{C}_{z_1})}$. Applying Theorem \ref{linearizzazione sul fibrato tangente}, it follows that $f_\alpha$ and $f_\beta$ are locally topologically conjugated in a neighborhood of $\mathbb{C}_{z_1}\subseteq\mathbb{C}^n$. Since $f_\alpha$ coincides with $f$ in a neighborhood of $0\in\mathbb{C}^n$, and $f_\beta=df_0$, we can conclude that $f$ and $df_0$ are locally topologically conjugated in a neighborhood of $0\in\mathbb{C}^n$. 

In case $(ii)$, up to a polynomial conjugacy, we have that, if $U\subseteq\mathbb{C}_{z_1}$ is a small neighborhood of the origin, then $f|_U$ is a $C^r$ function, with Taylor expansion at the origin $f|_U(z_1)=\lambda_1 z_1+\lambda_1 z_1^{kq+1}+o(|z_1|^{kq+1})$, with $k\geq1$ and $r>kq+1$. If $\phi:\{z_1\in\mathbb{C}_{z_1}\,\,\big|\,\,|z_1|<\rho\}\rightarrow\mathbb{C}_{z_1}$ is the holomorphic function introduced in Lemma \ref{traslazione all'infinito}, let us consider the holomorphic function $\hat{\phi}(z_1,\ldots,z_n):=(\phi(z_1), \lambda_2 z_2,\ldots,\lambda_n z_n)$, defined in a neighborhood of $0\in\mathbb{C}^n$. If $L:=df_0$, let us consider the constants $\gamma=\gamma(L)$, $\epsilon=\epsilon(L,\gamma)$ introduced in Propositions \ref{costruzione della varieta normalmente iperbolica}, \ref{coniugio dei differenziali}. Choosing $0<\eta<\rho$, let us define $\tilde{f}_\alpha(z):=\rho_{\frac{\eta}{2}}(z)f(z)+(1-\rho_{\frac{\eta}{2}}(z))\hat{\phi}(z)$, and $\tilde{f}_\beta(z):=\rho_{\frac{\eta}{2}}(z)\hat{g}(z)+(1-\rho_{\frac{\eta}{2}}(z))\hat{\phi}(z)$, where $\hat{g}(z):=(\lambda_1 (z_1+z_1^{kq+1}),\lambda_2 z_2,\ldots,\lambda_n z_n)$. These are $C^r$ functions, coincident for $\frac{\eta}{2}<|z|<\rho$; moreover, by Proposition \ref{risultato similar camacho}, there exists a homeomorphism $\Gamma:\mathbb{C}_{z_1}\rightarrow\mathbb{C}_{z_1}$, such that $\Gamma(z_1)=z_1$ for $|z_1|\geq\frac{\eta}{2}$, and $\Gamma\circ\tilde{f}_\alpha|_{\{z_1\in\mathbb{C}_{z_1}\,\,\big|\,\,|z_1|<\rho\}}=\tilde{f}_\beta|_{\{z_1\in\mathbb{C}_{z_1}\,\,\big|\,\,|z_1|<\rho\}}\circ\Gamma$. Let us define then $f_i(z):=\rho_\eta(z)\tilde{f}_i(z)+(1-\rho_\eta(z))L(z)$, for $i=\alpha,\beta$. Up to choosing $\eta$ suitably small, we can assume that $f_i:\mathbb{C}^n\rightarrow\mathbb{C}^n$ is a $C^r$ diffeomorphism, such that $\|f_i-L\|_{1,\mathbb{C}^n}<\epsilon$. Moreover, $\mathbb{C}_{z_1}\subseteq\mathbb{C}^n$ is $f_i$-invariant. By Proposition  \ref{costruzione della varieta normalmente iperbolica} it is normally hyperbolic for $f_i$, with splitting $T_{\mathbb{C}_{z_1}}\mathbb{C}^n=T\mathbb{C}_{z_1}\oplus(E^i)^s(\mathbb{C}_{z_1})\oplus(E^i)^u(\mathbb{C}_{z_1})$. By construction, we have also that $\Gamma\circ f_\alpha|_{\mathbb{C}_{z_1}}=f_\beta|_{\mathbb{C}_{z_1}}\circ\Gamma$. 
We are in the hypotheses of Proposition \ref{coniugio dei differenziali}, and there exists a homeomorphism $h:$ $(E^\alpha)^s(\mathbb{C}_{z_1})$ $\oplus$ $(E^\alpha)^u(\mathbb{C}_{z_1})$ $\rightarrow$ $(E^\beta)^s(\mathbb{C}_{z_1})$ $\oplus $ $(E^\beta)^u(\mathbb{C}_{z_1})$, which conjugates the restrictions $df_\alpha|_{(E^\alpha)^s(\mathbb{C}_{z_1})\oplus(E^\alpha)^u(\mathbb{C}_{z_1})}$ and $df_\beta|_{(E^\beta)^s(\mathbb{C}_{z_1})\oplus(E^\beta)^u(\mathbb{C}_{z_1})}$. Proceeding as in case $(i)$, we can  conclude that $f$ and $\hat{g}$ are locally topologically conjugated in a neighborhood of the origin; hence, up to a polynomial conjugacy, that $f$ and $g:z\mapsto(\lambda_1 z_1+z_1^{kq+1},\lambda_2 z_2,\ldots,\lambda_n z_n)$ are.
$\diamondsuit$


\end{document}